\documentclass{article}
\usepackage{amssymb}

\newtheorem{Theorem}{\sc Theorem}[section]
\newtheorem{Proposition}[Theorem]{\sc Proposition}
\newtheorem{Corollary}[Theorem]{\sc Corollary}
\newtheorem{Definition}[Theorem]{\sc Definition}
\newtheorem{Lemma}[Theorem]{\sc Lemma}
\newtheorem{Remark}[Theorem]{\sc Remark}
\input{tcilatex}
\input tcilatex
\begin{document}

\date{MARIA JOI\c{T}A}
\title{ON\ HILBERT\ MODULES\ OVER\ LOCALLY\ $C^{*}$-ALGEBRAS II}
\author{}
\maketitle

\begin{abstract}
$\!\!\!\!\!\!\!\!\!\!\!\!\!\!\!\!$ \ \ 

\ In this paper we study the unitary equivalence between Hilbert modules
over a locally $C^{*}$-algebra. Also, we prove a stabilization theorem for
countably generated modules over an arbitrary locally $C^{*}$-algebra and
show that a Hilbert module over a Fr\'{e}chet locally $C^{*}$-algebra is
countably generated if and only if the locally $C^{*}$-algebra of all
''compact'' operators has an approximate unit.


\end{abstract}

\section{Introduction}

The notion of Hilbert module over a locally $C^{*}$-algebra (inverse limit
of $C^{*}$-algebras) generalizes the notion of Hilbert $C^{*}$ -module by
allowing the inner product to take values in a locally $C^{*}$-algebra.

In [14], Phillips showed that most basic properties of Hilbert $C^{*}$%
-modules are valid for Hilbert modules over locally $C^{*}$-algebras, such
as a stabilization theorem for countably generated Hilbert modules over a
locally $C^{*}$-algebra whose topology is determined by a countably family
of $C^{*}$ -seminorms. In [5], it is proved a stabilization theorem for
countably bounded generated Hilbert modules over an arbitrary locally $C^{*}$%
-algebra. In this paper, we show that this theorem is true for countably
generated Hilbert modules over an arbitrary locally $C^{*}$-algebra (Theorem
4.2).

The paper is organized as follows. Section 2 contains some notation and
preliminaries. In Section 3 we present the necessary and sufficient
conditions for that two Hilbert $A$-modules to be isomorphic as Hilbert $A$%
-modules. Thus, we show that the Hilbert $A$-modules $E$ and $F$ are
unitarily equivalent if and only if the Hilbert $C^{*}$-modules $b(E)$ and $%
b(F)$ are unitarily equivalent (Corollary 3.7). Also, we prove that the
Hilbert $A$-modules $E$ and $F$ are unitarily equivalent if and only if
there is an adjointable operator $T$ from $E$ into $F$ such that $T$ and $%
T^{*}$ both have dense range (Proposition 3.8). In [1], Frank showed that
for a Banach $C^{*}$ -module $E$ over a $C^{*}$ -algebra carrying two $A$
-valued inner-products $\left\langle \cdot ,\cdot \right\rangle _1$ and $%
\left\langle \cdot ,\cdot \right\rangle _2$ which induce equivalent to the
given one norms on $E$ the appropriate $C^{*}$ -algebras $L_A(E,\left\langle
\cdot ,\cdot \right\rangle _1)$ and $L_A(E,\left\langle \cdot ,\cdot
\right\rangle _2)\;$are isomorphic if and only if there is a surjective
linear map $S$ from $(E,\left\langle \cdot ,\cdot \right\rangle _1)$ onto $%
(E,\left\langle \cdot ,\cdot \right\rangle _2)\;$such that $\left\langle
S\xi ,S\xi \right\rangle _2=\left\langle \xi ,\xi \right\rangle _1$ for all $%
\xi $ in $E.$ We extend this result in the context of Hilbert modules over
locally $C^{*}$ -algebras (Proposition 3.9). In Section 4, we show that the
stabilization theorem is also valid for countably generated Hilbert modules
over an arbitrary locally $C^{*}$-algebra. Using this theorem we show that
if $A$ is unital and $b(\mathcal{H}_A)=\mathcal{H}_{b(A)}$, then a Hilbert $%
A $-module $E$ is countably generated if and only if the Hilbert $b(A)$%
-module $b(E)$ is countably generated (Corollary 4.4). We know that a
Hilbert $C^{*}$-module $E$ is countably generated if and only if the $C^{*}$%
-algebra of all ''compact'' operators on $E$ has an approximate unit.
Finally we show that this result is valid in the case of Hilbert modules
over Fr\'{e}chet locally $C^{*}$-algebras (Proposition 4.5 and Corollary
4.5).

\section{Preliminaries}

\smallskip First we recall some results about locally $C^{*}$-algebras from
[14] and [2].

A locally $C^{*}$-algebra is a complete Hausdorff complex topological $*$%
-algebra $A$ whose topology is determined by its continuous $C^{*}$%
-seminorms in the sense that the net $\left\{ a_i\right\} _i$ converges to $%
0 $ if and only if the net $\left\{ p(a_i)\right\} _i\;$converges to $0$ for
every continuous $C^{*}$-seminorm $p$ on $A.$

A Fr\'{e}chet locally $C^{*}$-algebra is a locally $C^{*}$-algebra whose
topology is determined by a countable family of $C^{*}$ -seminorms.

Let $A$ be a locally $C^{*}$-algebra and let $S(A)$ be the set of all
continuous $C^{*}$-seminorms on $A.$ For $p\in S(A),$ $A_p=A/\ker p$ is a $%
C^{*}$-algebra in the norm induced by $p$, and for $p,q\in S(A),$ $p\geq q$
there is a canonical morphism $\pi _{pq}$ from $A_p$ onto $A_q$\ such that $%
\pi _{pq}(a+\ker p)=a+\ker p,$ $a\in A.$ Then $\{A_p,$ $\pi _{pq}\}_{p\geq
q,p,q\in S(A)}$ is an inverse system of $C^{*}$-algebras and the locally $%
C^{*}$-algebras $A$ and $\lim\limits_{\stackunder{p}{\leftarrow }}A_p$ are
isomorphic. The canonical map from $A$ onto $A_p$ will be denoted by $\pi
_p. $

The set $b(A)=\{a\in A;\left\| a\right\| _\infty =\sup \{p(a);p\in
S(A)\}<\infty \}$ of all bounded elements of $A$ with the $C^{*}$-norm $%
\left\| \cdot \right\| _\infty \;$is a $C^{*}$-algebra which is dense in $A.$

An approximate unit of \ $A$ is an increasing net $\{e_i\}_{i\in I}$ of
positive elements of $A\;$such that: $p(e_i)\leq 1$ for all $i\in I\;$and
for all $p\in S(A)$; $p(ae_i-a)\rightarrow 0$ and $p(e_ia-a)\rightarrow 0$
for all $p\in S(A)$ and for all $a\in A.$ Any locally $C^{*}$-algebra has an
approximate unit.

Now we will recall some results about Hilbert modules over locally $C^{*}$%
-algebras from [14], [5], [6] and [9].

\begin{Definition}
A pre-Hilbert $A$-module is a complex vector space\ $E$ which is also a
right $A$-module, compatible with the complex algebra structure, equipped
with an $A$-valued inner product $\left\langle \cdot ,\cdot \right\rangle
:E\times E\rightarrow A\;$which is $\Bbb{C}$- and $A$-linear in its second
variable and satisfies the following relations:

\ \ \ \ \ \ \ $\;\;(i)\;\;\left\langle \xi ,\eta \right\rangle
^{*}=\left\langle \eta ,\xi \right\rangle \;\;$for every $\xi ,\eta \in E;$

\ \ \ \ \ \ \ \ \ $(ii)\;\;\left\langle \xi ,\xi \right\rangle \geq 0\;\;$%
for every $\xi \in E;$

\ \ \ \ \ \ \ \ \ $(iii)\;\left\langle \xi ,\xi \right\rangle =0\;\;$if and
only if $\xi =0.$

\ We say that $E\;$is a Hilbert $A$-module if $E\;$is complete with respect
to the topology determined by the family of seminorms $\{\overline{p}%
_E\}_{p\in S(A)}\;$where $\overline{p}_E\left( \xi \right) =\sqrt{p\left(
\left\langle \xi ,\xi \right\rangle \right) },\xi \in E.$
\end{Definition}

\begin{Remark}
If $A\;$is not unital and $E$\ is a Hilbert $A$-module, then $E$ becomes a
Hilbert $A^{+}$-module, where $A^{+}$ is the unitization of $A$, if we
define $\xi 1=\xi .$ Moreover, $EA$ is a pre-Hilbert $A$ -module which is
dense in $E.$
\end{Remark}

Let $E\;$be a Hilbert $A$-module.\ For $p\in S(A),\;\mathcal{E}_p=\{\xi \in
E;p(\left\langle \xi ,\xi \right\rangle )=0\}\;$is a closed submodule of $%
E\; $and $E_p=$ $E/\mathcal{E}_p\;$is a Hilbert $A_p$-module with: $(\xi +%
\mathcal{E}_p)\pi _p(a)=\xi a+\mathcal{E}_p\;$and $\left\langle \xi +%
\mathcal{E}_p,\eta +\mathcal{E}_p\right\rangle =\pi _p\left( \left\langle
\xi ,\eta \right\rangle \right) $. The canonical map from $E\;$onto $E_p$
will be denoted by $\sigma _p^E.$

For $p,q\in S(A)$ with $p\geq q,$ there is a canonical morphism of vector
spaces $\sigma _{pq}^E\;$from $E_p\;$onto $E_q\;$such that $\sigma
_{pq}^E\left( \sigma _p^E\left( \xi \right) \right) =\sigma _q^E\left( \xi
\right) ,\;\xi \in E.\;$Then $\{E_p;$\ $A_p;$ $\sigma _{pq}^E\}_{p\geq
q;p,q\in S(A)}$ is an inverse system of Hilbert $C^{*}$-modules in the
following sense: $\sigma _{pq}^E\left( \sigma _p^E\left( \xi \right) \pi
_p\left( a\right) \right) =\sigma _{pq}^E\left( \sigma _p^E\left( \xi
\right) \right) \pi _{pq}\left( \pi _p\left( a\right) \right) ,\xi \in
E,a\in A;$ $\left\langle \sigma _{pq}^E\left( \sigma _p^E\left( \xi \right)
\right) ,\sigma _{pq}^E\left( \sigma _p^E\left( \eta \right) \right)
\right\rangle =\pi _{pq}\left( \pi _p\left( \left\langle \xi ,\eta
\right\rangle \right) \right) ,\xi ,\eta \in E;\sigma _{pp}^E=$id$_{E_p};$ $%
\sigma _{qr}^E\circ \sigma _{pq}^E=\sigma _{pr}^E\;$if $p\geq q\geq r,$ and $%
\lim\limits_{\stackunder{p}{\leftarrow }}E_p$ is a Hilbert $A$-module with: $%
\left( \sigma _p^E\left( \xi \right) \right) _p\left( \pi _p\left( a\right)
\right) _p$ $=\left( \sigma _p^E\left( \xi a\right) \right) _p$ and $%
\left\langle \left( \sigma _p^E\left( \xi \right) \right) _p,\left( \sigma
_p^E\left( \eta \right) \right) _p\right\rangle =\left( \left\langle \sigma
_p^E\left( \xi \right) ,\sigma _p^E\left( \eta \right) \right\rangle \right)
_p$.\ Moreover, the Hilbert $A$-modules\ $E$ and $\lim\limits_{\stackunder{p%
}{\leftarrow }}E_p$\ may be identified.

The set $b(E)=\{\xi \in E;\left\| \xi \right\| _\infty =\sup \{\overline{p}%
(\xi );p\in S(A)\}<\infty \}$ of all bounded elements of $E$ is a Hilbert $%
b(A)$-module.

Let $E\;$and $F$\ be Hilbert $A$-modules. We will denote by $L_A(E,F)$ the
set of all adjointable $A$-module homomorphisms from $E$ into $F,\ $and we
write $L_A(E)\;$for $L_A(E,E).$

For $p\in S(A)$, define $\left( \pi _p\right) _{*}\;$from $L_A(E,F)$ into $%
L_{A_p}(E_p,F_p)$\ by $\left( \pi _p\right) _{*}\left( T\right) \;(\xi +%
\mathcal{E}_p)=T\xi +\mathcal{F}_p,$ $\xi \in E$ and for $p,q\in S(A)$ with $%
p\geq q$, define $\left( \pi _{pq}\right) _{*}$ from $L_{A_p}(E_p,F_p)$ into 
$L_{A_q}(E_q,F_q)$ by $\left( \pi _{pq}\right) _{*}\left( T_p\right) \left(
\xi +\mathcal{E}_q\right) =\sigma _{pq}^F\left( T_p\left( \xi +\mathcal{E}%
_p\right) \right) ,\xi \in E.\;$Then $\left\{ L_{A_p}(E_p,F_p);\left( \pi
_{pq}\right) _{*}\right\} _{p\geq q,p,q\in S(A)\;}$is an inverse system of
Banach spaces and if we consider on $L_A(E,F) $ the topology determined by
the family of seminorms $\left\{ \widetilde{p}\right\} _{p\in S(A)},$ where $%
\widetilde{p}(T)=\left\| (\pi _p)_{*}(T)\right\| _{L_{A_p}(E_p,F_p)}$ , $%
\left\| \cdot \right\| _{L_{A_p}(E_p,F_p)}$ being the operator norm on $%
L_{A_p}(E_p,F_p),$ $L_A(E,F) $ may be identified to $\lim\limits_{%
\stackunder{p}{\leftarrow }}L_{A_p}(E_p,F_p).$ Thus topologized, $L_A(E)$
becomes a locally $C^{*}$-algebra.

The set $b(L_A(E,F))=\{T\in L_A(E,F);\left\| T\right\| _\infty =\sup \{%
\widetilde{p}(T);p\in S(A)\}<\infty \}$ of all bounded elements in $L_A(E,F)$
is a Banach space with respect to the norm $\left\| \cdot \right\| _\infty .$

By definition, the set of all ''compact'' operators $K_A(E)$ on $E$ is
defined as the closure of the set of all finite linear combinations of the
operators 
\[
\left\{ \theta _{\xi ,\eta };\theta _{\xi ,\eta }\left( \zeta \right) =\xi
\left\langle \eta ,\zeta \right\rangle ,\xi ,\eta ,\zeta \in E\right\} . 
\]
It is a locally $C^{*}$-subalgebra and a two-sided ideal of $L_A(E)$ and
moreover $K_A(E)$ may be identified to $\lim\limits_{\stackunder{p}{%
\leftarrow }}K_{A_p}(E_p).$

\ If $E$ and $F$ are Hilbert $A$ -modules, then we can form the direct sum $%
E\oplus F.$ This is a Hilbert $A$ -module with $(\xi \oplus \eta )a=\xi
a\oplus \eta a,$ $\xi \in E,\eta $ $\in F,$ $a\in A$ and $\left\langle \xi
_1\oplus \eta _1,\xi _2\oplus \eta _2\right\rangle =\left\langle \xi _1,\xi
_2\right\rangle +\left\langle \eta _1,\eta _2\right\rangle ,$ $\xi _1,\xi
_2\in E,\eta _1,\eta _2\in F.$

\section{\protect\smallskip Unitarily equivalent Hilbert $A$-modules}

Let $A$ be a locally $C^{*}$ -algebra and let\ $E$\ and\ $F$\ be Hilbert $A$%
\ -modules.

\begin{Definition}
An operator $U\in L_A(E,F)$\ is said to be unitary if $U^{*}U=$id$_E$\ and $%
UU^{*}=$id$_F.$
\end{Definition}

\begin{Remark}
If\ $U$ is an element in\ $L_A(E,F),$\ then $U$\ is unitary if and only if $%
(\pi _p)_{*}(U)$\ is a unitary element in $L_{A_p}(E_p,F_p)$ for each $p\in
S(A)$.
\end{Remark}

From [12], Theorem 3.5 and Remark 3.2, we obtain:

\begin{Proposition}
Let $E$ and $F$ be Hilbert $A$ -modules and let $U$ be a linear map from $E$
into $F.$ Then the following statements are equivalent:

\begin{enumerate}
\item  $U$ is unitary;

\item  $\left\langle U\xi ,U\xi \right\rangle =\left\langle \xi ,\xi
\right\rangle $ for all $\xi \in E$ and $U$ is surjective;

\item  $\overline{p}_F\left( U\xi \right) =\overline{p}_E\left( \xi \right) $
for all $\xi \in E$ and $U$ is a surjective $A$ -linear map.
\end{enumerate}
\end{Proposition}

\begin{Definition}
We say that the Hilbert $A$-modules $E$ and $F$ are unitarily equivalent,
and we write $E\thickapprox F$, if there is a unitary operator in $L_A(E,F)$.
\end{Definition}

\begin{Remark}
If $E\thickapprox F,$ then $E_p\thickapprox F_p$ for all $p\in S(A)$.
\end{Remark}

\begin{Proposition}
The set $\mathcal{U}_A(E,F)$ of all unitary operators from $E$ to $F$ is
isomorphic as set with $\mathcal{U}_{b(A)}(b(E),b(F))$, the set of all
unitary operators from $b(E)$ to $b(F)$.
\end{Proposition}

\proof%
By [9] Theorem 3.9, the map $\Psi :b(L_A(E,F))\rightarrow
L_{b(A)}(b(E),b(F)) $ defined by $\Psi (T)=T|_{b(E)}$ is an isometrically
isomorphism of Banach spaces. Clearly $\mathcal{U}_A(E,F)\subseteq
b(L_A(E,F)).$ It is not hard to check that the restriction of $\Psi $ on $%
\mathcal{U}_A(E,F)$ is an isomorphism of set from $\mathcal{U}_A(E,F)$ onto $%
\mathcal{U}_{b(A)}(b(E),b(F))$.%
\endproof%

\begin{Corollary}
Let $E$ and $F$ be Hilbert $A$ -modules. Then $E$ and $F$ are unitarily
equivalent if and only if $b(E)$ and $b(F)$ are unitarily equivalent.
\end{Corollary}

The following proposition is a generalization of Proposition 3.8, [12] in
the context of Hilbert modules over locally $C^{*}$-algebras.

\begin{Proposition}
The Hilbert $A$-modules $E$ and $F$ are unitarily equivalent if and only if
there is an element $T$ in $L_A(E,F)\;$such that $T\;$and $T^{*}\;$have
dense range.
\end{Proposition}

\proof%
If $E\thickapprox $ $F$, then there is a unitary operator $U$ in $L_A(E,F)$.
Since $U$ is unitary, $U$ and $U^{*}$ are surjective and so $U$ and $U^{*}$
have dense range.

Conversely, if $T$ and $T^{*}$ have dense range, then, by [8], Theorem 2.8, $%
T$ has a polar decomposition. Therefore $T=U|T|,$ where $U$ is a partial
isometry in $L_A(E,F)$ such that $U^{*}U$ is the projection of $E$ on $%
\overline{TE}$ and $UU^{*}$ is the projection of $F$ on $\overline{|T|E}.$
Since $T$ and $T^{*}$ have dense range and $\overline{T^{*}E}=\overline{|T|E|%
},$ $U$ is a unitary operator in $L_A(E,F)$ and so $E\thickapprox F$.%
\endproof%

Let $E$ be a complex vector space which is also right $A$ -module,
compatible with the structure of complex algebra and equipped with two an $A$
-valued inner-products $\left\langle \cdot ,\cdot \right\rangle _1$ and $%
\left\langle \cdot ,\cdot \right\rangle _2$ which induce either a structure
of Hilbert $A$ -module on $E.$ We denote by $E_1$ the Hilbert $A$ -module $%
\left( E,\left\langle \cdot ,\cdot \right\rangle _1\right) $ and by $E_2$
the Hilbert $A$ -module $\left( E,\left\langle \cdot ,\cdot \right\rangle
_2\right) .$

The following proposition is a generalization in the context of Hilbert
modules over locally $C^{*}$ -algebras of a result of M. Frank [1].

\begin{Proposition}
Let $E$ be as above. Then the following statements are equivalent:

\begin{enumerate}
\item  The Hilbert $A$ -modules $E_1$ and $E_2$ are unitarily equivalent.

\item  The locally $C^{*}$ -algebras $K_A(E_1)$ and $K_A(E_2)$ are
isomorphic.

\item  The locally $C^{*}$ -algebras $L_A(E_1)$ and $L_A(E_2)$ are
isomorphic.

\item  The $C^{*}$ -algebras $L_{b(A)}(b(E_1))$ and $L_{b(A)}(b(E_2))$ are
isomorphic.

\item  The $C^{*}$ -algebras $K_{b(A)}(b(E_1))$ and $K_{b(A)}(b(E_2))$ are
isomorphic.

\item  The Hilbert $b(A)$ -modules $b(E_1)$ and $b(E_2)$ are isometrically
isomorphic as Banach $b(A)$ -modules.

\item  The Hilbert $b(A)$ -modules $b(E_1)$ and $b(E_2)$ are unitarily
equivalent.
\end{enumerate}
\end{Proposition}

\proof%
$1.\Rightarrow 2.$ Since $E_1$ and $E_2$ are unitarily equivalent, there is
a unitary operator $U$ in $L_A(E_1,E_2)$. It is not hard to check that the
map $\Phi $ from $K_A(E_1)$ to $K_A(E_2)$ defined by $\Phi (T)=UTU^{*}$ is
an isomorphism of locally $C^{*}$ -algebras.

$2.\Rightarrow 3.$ Let $\Phi $ be an isomorphism of locally $C^{*}$
-algebras from $K_A(E_1)$ onto $K_A(E_2).$ By [4], Lemmas 2.4, 2.7 and
Corollary 2.6, there is a unique isomorphism of locally $C^{*\text{ }}$%
-algebras $\overline{\Phi }:M(K_A(E_1))\rightarrow $ $M(K_A(E_2)),$ where $%
M(K_A(E_i))$ denotes the locally $C^{*}$ -algebra of all multipliers of $%
K_A(E_i),$ $i=1,2,$ such that $\overline{\Phi }|_{K_A(E_1)}=\Phi .$

On the other hand, the locally $C^{*}$ -algebras $M(K_A(E_1))$ and $L_A(E_1)$
are isomorphic as well as $M(K_A(E_2))$ and $L_A(E_2)$ ( [14],Theorem 4.2).
Therefore the locally $C^{*}$ -algebras $L_A(E_1)$ and $L_A(E_2)$ are
isomorphic.

$3.\Rightarrow 4.$ If the locally $C^{*}$ -algebras $L_A(E_1)$ and $L_A(E_2)$
are isomorphic, then the $C^{*}$ -algebras $b(L_A(E_1))$ and $b(L_A(E_2))$
are isomorphic ([14], Corollary 1.13). But $b(L_A(E_i))$ is isomorphic with $%
L_{b(A)}(b(E_i)),$ $i=1,2$ ( [7], Theorem 3.3 ). Therefore the $C^{*}$
-algebras $L_{b(A)}(b(E_1))$ and $L_{b(A)}(b(E_2))$ are isomorphic.

The implications $4.\Rightarrow 5.\Rightarrow 6.$ were proved in [1], the
equivalence $6.\Leftrightarrow 7.$ was proved in [11] and the implication $%
7.\Rightarrow 1.$ was showed in Corollary 3.7.%
\endproof%

As in the case of Hilbert $C^{*}$ -modules, the Hilbert $A$ -module $%
\mathcal{H}_A=\{(a_n)_n;\sum\limits_{n=1}^\infty a_n^{*}a_n$ is convergent
in $A\}$ plays a special role in the theory of Hilbert modules over locally $%
C^{*}$ -algebras. For each $p\in S(A),$ the Hilbert $A_p$ -modules $\mathcal{%
H}_{A_p}$ and $(\mathcal{H}_A)_p$ are unitarily equivalent, and moreover,
the Hilbert $A$ -modules $\mathcal{H}_A$ and $\lim\limits_{\stackunder{p}{%
\leftarrow }}\mathcal{H}_{A_p}$ unitarily equivalent (see [14], Section 4).

\begin{Lemma}
Let $A$ be a non unital locally $C^{*}$ -algebra and let $E$ be a Hilbert $A$
-module. If $A^{+}$ is the unitization of $A,$ then the Hilbert $A$ -modules 
$\mathcal{H}_A$ and $\overline{\mathcal{H}_{A^{+}}A}$ are unitarily
equivalent as well as $E\oplus \mathcal{H}_A$ and $\ \overline{\left(
E\oplus \mathcal{H}_{A^{+}}\right) A}.$
\end{Lemma}

\proof%
It is not hard to check that the map $U$ from $\mathcal{H}_{A^{+}}A$ to $%
\mathcal{H}_A$ defined by $U((a_n)_nb)=(a_nb)_n$ extends by continuity to a
unitary from $\overline{\mathcal{H}_{A^{+}}A}$ to $\mathcal{H}_A$ and the
map $V$ from $\left( E\oplus \mathcal{H}_{A^{+}}\right) A$ to $E\oplus 
\mathcal{H}_A$ defined by $V((\xi \oplus (a_n)_n)b)=\xi \oplus (a_nb)_n$
extends to a unitary from $\overline{\left( E\oplus \mathcal{H}%
_{A^{+}}\right) A}$ to $E\oplus \mathcal{H}_A.$%
\endproof%

\begin{Remark}
Let $A$ be a locally $C^{*}$ -algebra. Then $\mathcal{H}_{b(A)}$ is a closed
submodule of $b(\mathcal{H}_A)=\{(a_n)_n;\dsum\limits_na_n^{*}a_n$ converges
in $A$ to an element in $b(A)\}.$ In general $\mathcal{H}_{b(A)}$ does not
coincides with $b(\mathcal{H}_A).$

\textsc{Example} $1$. Let $A=C(\Bbb{Z}^{+}).$ Then $A$ equipped with the
topology determined by the family of $C^{*}$ -seminorm $\{p_n\}_n,$ where $%
p_n(f)=\sup \{|f(k)|;1\leq k\leq n\}$ is a locally $C^{*}$ -algebra.

For each positive integer $n,$ we consider the function $f_n$ from $\Bbb{Z}%
^{+}$ to $\Bbb{C}$ defined by 
\[
f_n(m)=\left\{ 
\begin{array}{ll}
1 & \text{if }m=n \\ 
0 & \text{if }m\neq n
\end{array}
\right. . 
\]
It is easy to check that $\tsum\limits_n\left| f_n\right| ^2$ is convergent
in $A.$ Hence $(f_n)_n$ is an element in $b(H_A).$

Since 
\[
\sup \{\sum_{n\geq n_0}\left| f_n\right| ^2(m);m\in \Bbb{Z}^{+}\}=1 
\]
for any positive integer $n_0,$ $\tsum\limits_n\left| f_n\right| ^2$ is not
convergent in $b(A),$ and so $(f_n)_n$ $\notin H_{b(A)}.$ Therefore $%
H_{b(A)}\varsubsetneqq b(H_A).$

\textsc{Example} $2$. Let $A=C_{cc}[0,1]$ be the set of all $\func{complex}$
continuous functions on $[0,1]$ endowed with the topology ''cc'' of uniform
convergence on the countable compact subsets of $[0,1].$ Then $A$ is a
locally $C^{*}$ -algebra .

If $(f_n)_n$ $\in b(H_A),$ then $\tsum\limits_n\left| f_n\right| ^2$ is
convergent in $A,$ and by Dini's theorem it is uniformly convergent.
Therefore $\tsum\limits_n\left| f_n\right| ^2$ is convergent in $b(A)$ and
so $(f_n)_n$ $\in H_{b(A)}.$ Hence $H_{b(A)}=b(H_A).$
\end{Remark}

\section{Countably generated Hilbert $A$-modules}

Let $A$ be a locally $C^{*}$ -algebra and let $E$ be a Hilbert $A$ -module.
A subset $G$ of $E$ is a generating set for $E$ if the closed submodule of $%
E $ generated by $G$ is the whole of $E.$ We say that $E$ is countably
generated if it has a countable generating set.

\begin{Lemma}
If $E$ is countably generated then it has a generating set contained in $%
b(E).$ Moreover, $E_p$ is countably generated for each $p\in S(A)$ .
\end{Lemma}

\proof%
Let $\{\xi _n;$ $n=1,2,...\}$ be a generating set for $E.$ According to [3],
Proposition 3.2, for each positive integer $n$ there is a sequence $\{\xi
_n^m\}_m$ in $b(E)$ such that $\xi _n^m\rightarrow \xi _n.$ Then $\{\xi
_n^m;n,m=1,2,...\}$ is a generating set for $E.$

Let $p\in S(A).$ Since the canonical map $\sigma _p$ from $E$ onto $E_p$ is
a surjective continuous map and $\sigma _p(\xi a)=\sigma _p(\xi )\pi _p(a)$
for all $\xi \in E$ and $a\in A,$ $\{\sigma _p(\xi _n);n=1,2,...\}$ is a
generating set for $E_p,$ and so $E_p$ is countably generated.%
\endproof%

Now, using Lemma 4.1 and Kasparov's theorem for countably generated Hilbert $%
C^{*}$-modules, we prove a stabilization theorem for countably generated
Hilbert modules over locally $C^{*}$-algebras.

\begin{Theorem}
If $A$ is a locally $C^{*}$ -algebra and $E$ is a countably generated
Hilbert $A$ -module then $\mathcal{H}_A\thickapprox E\oplus \mathcal{H}_A.$
\end{Theorem}

\proof%
First we suppose that $A$\ has a unity $1.$

For each positive integer $n$\ we denote by $e_n$\ the element in $\mathcal{H%
}_A$\ which has all the components zero except at the $n^{\text{th}}$\
component which is $1.$

Let $\{\xi _n;n=1,2,...$\ $\}$\ be a generating set for $E$\ with each
element repeated infinitely often. According to Lemma 4.1, we can suppose
that $\xi _n\;$is bounded and $\left\| \xi _n\right\| _\infty \leq 1$\ for
all positive integer $n$.\ Then $\{\sigma _p(\xi _n);n=1,2,...\}$ is a
generating set for $E_p$ with each element repeated infinitely often and $%
\left\| \sigma _p(\xi _n)\right\| _{\overline{p}}\leq 1$\ for all positive
integer $n,$ and for each positive integer $n$, $\sigma _p^{\mathcal{H}%
_A}(e_n)$ is the element in $\mathcal{H}_{A_p}$ which has all the components
zero except at the $n^{\text{th}}$\ component which is $1.$ By Kasparov's
theorem for countably generated Hilbert $C^{*}$ -modules, there is a unitary
operator $U_p$ from $\mathcal{H}_{A_p}$ to $E_p\oplus \mathcal{H}_{A_p}.$
Moreover, according to the proof of this theorem (see, for instance [],
Theorem 6.2), we can choose $U_p$ such that $U_p|T_p|=T_p,$ where $T_p$ is
an element in $L(\mathcal{H}_{A_p,}E_p\oplus \mathcal{H}_{A_p})$ defined by 
\[
T_p\mathbf{=}\tsum\limits_n\theta _{2^{-n}\sigma _p^E(\xi _n)\oplus
4^{-n}\sigma _p^{\mathcal{H}_A}(e_n),\sigma _p^{\mathcal{H}_A}(e_n)} 
\]
and 
\[
\left| T_p\right| =(T_p^{*}T_p)^{\frac 12}. 
\]
It is not hard to check that $(T_p)_p\;$is a coherent sequence in $L(%
\mathcal{H}_{A_p},E_p\oplus \mathcal{H}_{A_p})$\ and so $\left( \left|
T_p\right| \right) _p$ is a coherent sequence in $L(\mathcal{H}_{A_p}).$\ 

Knowing that for each $p\in S(A),$ $|T_p|$ has dense range, from 
\[
(\pi _{pq})_{*}(U_p)|T_q|=(\pi _{pq})_{*}(U_p|T_p|)=(\pi
_{pq})_{*}(T_p)=T_q=U_q|T_q| 
\]
for all $p,q\in S(A)$ with $p\geq q,$\ we conclude that $(U_p)_p$\ is a
coherent sequence in $L(\mathcal{H}_{A_p},E_p\oplus \mathcal{H}_{A_p}).$
Therefore the Hilbert $A$ -modules $\lim\limits_{\stackunder{p}{\leftarrow }}%
\mathcal{H}_{A_p}$ and $\lim\limits_{\stackunder{p}{\leftarrow }}(E_p\oplus 
\mathcal{H}_{A_p})$ are unitarily equivalent and since the Hilbert $A$
-modules $\mathcal{H}_A$ and $\lim\limits_{\stackunder{p}{\leftarrow }}%
\mathcal{H}_{A_p{}}$ are unitarily equivalent as well as $E\oplus \mathcal{H}%
_A$ and $\lim\limits_{\stackunder{p}{\leftarrow }}(E_p\oplus \mathcal{H}%
_{A_p}),$ the theorem is proof in this case.

If $A$ is not unital, let $A^{+}$ be the unitization of $A.$ Since $E$ can
be regarded as a Hilbert $A^{+}$ -module, by the first part of the proof
there is a unitary operator $U^{+}$ from $\mathcal{H}_{A^{+}}$ to $E\oplus 
\mathcal{H}_{A^{+}}.$ Let $U$ be the restriction of $U^{+}$ on $\overline{%
\mathcal{H}_{A^{+}}A}.$ Then $U$ is a unitary operator from $\overline{%
\mathcal{H}_{A^{+}}A}$ to $\overline{(E\oplus \mathcal{H}_{A^{+}})A}$ and by
Lemma 3.10, $\mathcal{H}_A$ is unitarily equivalent with $E\oplus \mathcal{H}%
_A$. 
\endproof%

\begin{Remark}
If $E$ is a Hilbert $A$ -module such that $b(E)$ is countably generated,
then clearly $E$ is countably generated.

If $E$ is a countable generated Hilbert $A$ -module, then is $b(E)$ a
countably generated Hilbert $b(A)$ -module? We solve this problem in a
particular case when $b(\mathcal{H}_A)=\mathcal{H}_{b(A)}.$
\end{Remark}

\begin{Corollary}
Let $A$ be a unital locally $C^{*}$-algebra such that $b(\mathcal{H}_A)=%
\mathcal{H}_{b(A)}$ and let $E$ be a Hilbert $A$ -module. If $E$ is
countably generated, then $b(E)$ is countably generated.
\end{Corollary}

\proof%
From Theorem 4.2 and Corollary 3.7, we conclude that $b\left( E\oplus 
\mathcal{H}_A\right) \thickapprox b\left( \mathcal{H}_A\right) .$ It is easy
to see that $b(E\oplus \mathcal{H}_A)\thickapprox b(E)\oplus b(\mathcal{H}%
_A) $ and so $b(E)\oplus b(\mathcal{H}_A)\thickapprox b(H_A).$ But $b(H_A)\ $%
is countably generated, since $A$ is unital and so $b(A)$ is unital, and
consequently $b(E)$ is countably generated.%
\endproof%

It is well-known that a Hilbert $C^{*}$ -module $E$ over a $C^{*}$ -algebra $%
A$ is countably generated if and only if $K_A(E)$ has a countably
approximate unit ([10]). We show that this result is valid for Hilbert
modules over Fr\'{e}chet locally $C^{*}$ -algebras.

\begin{Proposition}
Let $A$ be a locally $C^{*}$ -algebra and let $E$ be a Hilbert $A$ -module.
Then:

\begin{enumerate}
\item  If $E$ is countably generated, then $K_A(E)$ has a countable
approximate unit.

\item  If $A$ is a Fr\'{e}chet locally $C^{*}$ -algebra and if $K_A(E)$ has
a countable approximate unit, then $E$ is countably generated.
\end{enumerate}
\end{Proposition}

\proof%
$1.$ According to Remark 2.2 we can suppose that $A$ is unital.

Let $p\in S(A)$ and $T_p=\sum\limits_n2^{-n}\theta _{\sigma _p^{\mathcal{H}%
_A}(e_n),\sigma _p^{\mathcal{H}_A}(e_n)}$, where, for each positive integer $%
n$,$\ e_n$\ denotes the element in $\mathcal{H}_A$\ which has all the
components zero except at the $n^{\text{th}}$\ component which is $1.\;$%
Since $T_p$ is\ a strictly positive element in $K_{A_p}(\mathcal{H}_{A_p})$
(see [12], pp. 66), $\{T_p^{\frac 1n}\}_n$ is an approximate unit of $%
K_{A_p}(\mathcal{H}_{A_p})$ (see [13], 3.10.4-3.10.5). It is easy to check
that for each $n,$ $(T_p^{\frac 1n})_p$ is a coherent sequence in $K_{A_p}(%
\mathcal{H}_{A_p})$, and so $\{T^{\frac 1n}\}_n,$ where $T^{\frac
1n}=(T_p^{\frac 1n})_p$is an approximate unit for $K_A(\mathcal{H}_A).$

If $E$ is countably generated, by Theorem 4.2,$\;E$ may be identified with a
complemented submodule of $\mathcal{H}_A.$ Let $P\in L(\mathcal{H}_A,E)$ be
the projection from $\mathcal{H}_A$ onto $E.$ It is easy to verify that $%
\{PT^{\frac 1n}P^{*}\}_n$ is an approximate unit for $K_A(E).$

$2.$ Let $\{V_n\}_n$ be an approximate unit of $K_A(E)$ and let $%
T=\sum\limits_n2^{-n}$ $V_n.$ Clearly $T\in K_A(E).$ We will show that $T$
has dense range.

Let $p\in S(A).$ Since $\{(\pi _p)_{*}(V_n)\}_n$ is an approximate unit of $%
K_{A_p}(E_p),$ $(\pi _p)_{*}(T)$ has dense range (see [13], 3.10.4-3.10.5
and [3], Lemma 1.1.21). Thus we have

\begin{eqnarray*}
\overline{TE} &=&\lim\limits_{\stackunder{p}{\leftarrow }}\overline{\sigma
_p^E(TE)}=\lim\limits_{\stackunder{p}{\leftarrow }}\overline{(\pi
_p)_{*}(T)(E_p)} \\
&=&\lim\limits_{\stackunder{p}{\leftarrow }}E_p=E.
\end{eqnarray*}
where $\overline{X}$ means the closure of the vector subspace $X$ with
respect to the topology induced by the inner product. Therefore $T$ has
dense range.

Let $\{p_n\}_n$ be a cofinal subset of $S(A).$ Since $T$ $\in K_A(E)$, for
each positive integer $n,$ there are $\xi _1^n,...,\xi _{m_n}^n$ $\eta
_1^n,...,\eta _{m_n}^n$ in $E$ such that 
\[
\widetilde{p_n}\left( T-\sum\limits_{k=1}^{m_n}\theta _{\xi _k^n,\eta
_k^n}\right) <\frac 1{2n}. 
\]
We show that $\{\xi _k^n;1\leq k\leq m_n,n=1,2,...\}$ is a system of
generators for $E.$

Let $\xi $ $\in $ $E,$ $\varepsilon >0\;$and let $n_0$ be a positive
integer. Since $T$ has dense range, there is $\eta \in E$ such that $\;%
\overline{p_{n_0}}(\xi -T\eta )$ $<\frac \varepsilon 2.$ Let $n=\max \{n_0,[%
\frac{\overline{p_{n_0}}(\eta )}\varepsilon ]+1\},$ where $[t]$ means the
integer part of the positive number $t.$ Then $p_{n_0}\leq p_n$ and 
\begin{eqnarray*}
\overline{p_{n_0}}\left( \xi -\sum\limits_{k=1}^{m_n}\xi _k^n\left\langle
\eta _k^n,\eta \right\rangle \right) &\leq &\overline{p_{n_0}}\left( \xi
-T\eta \right) +\overline{p_{n_0}}\left( T\eta
-\sum\limits_{k=1}^{m_n}\theta _{\xi _k^n,\eta _k^n}(\eta )\right) \\
&<&\frac \varepsilon 2+\overline{p_{n_0}}\left( \eta \right) \widetilde{%
p_{n_0}}\left( T-\sum\limits_{k=1}^{m_n}\theta _{\xi _k^n,\eta _k^n}\right)
\\
&<&\frac \varepsilon 2+\overline{p_{n_0}}\left( \eta \right) \widetilde{p_n}%
\left( T-\sum\limits_{k=1}^{m_n}\theta _{\xi _k^n,\eta _k^n}\right)
<\varepsilon .
\end{eqnarray*}
This shows that $\{\xi _k^n;1\leq k\leq m_n,n=1,2,...\}$ is a system of
generators for $E$ and therefore $E$ is countably generated. 
\endproof%

\begin{Corollary}
Let $A$ be a Fr\'{e}chet locally $C^{*}$ -algebra and let $E$ be a Hilbert $%
A $ -module. Then $E$ is countably generated if and only if $K_A(E)$ has a
countable approximate unit.
\end{Corollary}

\smallskip Department of Mathematics, Faculty of Chemistry, University of
Bucharest, Bd. Regina Elisabeta nr.4-12, Bucharest, Romania, \smallskip
mjoita@fmi.unibuc.ro

\end{document}